\documentclass[12pt,intlimits,reqno]{article}

\usepackage{amssymb}
\usepackage{epsfig}
\usepackage{graphicx}
\usepackage{lscape}
\usepackage{amsmath}
\usepackage{amsfonts}
\usepackage{amscd}
\usepackage{cite}
\usepackage{amsthm}

\usepackage[all,cmtip]{xy}

\newtheorem{theorem}{Theorem}

\begin{document}

\thispagestyle{empty}

\begin{center}
    {\huge \textbf{Parametric Solutions for}} \\[0.4cm]
    {\huge \textbf{a Nearly-Perfect Cuboid}} \\[1.5cm]
    {\Large \textbf{Mamuka Meskhishvili}}
\end{center}

\vskip+1.5cm

\begin{center}
    {\fontsize{14}{16pt}\selectfont \textbf{Abstract} }
\end{center}
\vskip+0.5cm

We consider nearly-perfect cuboids (NPC), where the only irrational is one of the face diagonals. Obtained are three rational pa\-ra\-met\-ri\-za\-ti\-ons for NPC with one parameter.

\vskip+1.5cm

\textbf{Keywords.}
    Perfect cuboid, perfect box, nearly-perfect cuboid, semi-perfect cuboid, Euler brick, rational cuboid, rational pa\-ra\-met\-ri\-za\-ti\-on.

\vskip+0.7cm

\textbf{2010 AMS Classification.}
    11D41, 11D72, 14G05.

\newpage

\section{Introduction}
\vskip+0.5cm

For three centuries, no one has been able to prove the existence or otherwise of a perfect cuboid (PC) -- a rectangular parallelepiped having integer (rational): three sides, three face (surface) diagonals and a space (body) diagonal.

If among the seven:
\vskip+0.02cm
{\fontsize{14}{16pt}\selectfont
    $$  a,\;b,\;c,\;d_{ab},\;d_{bc},\;d_{ac},\;d_s     $$
 }
\vskip+0.2cm
\noindent
only one is irrational, this is called a nearly-perfect cuboid (NPC). The    \linebreak NPC are divided into three types:
\begin{itemize}
\item[1.] The space diagonal {\fontsize{14}{16pt}\selectfont $d_s$} is irrational,

\item[2.] One of the sides {\fontsize{14}{16pt}\selectfont $a$}, {\fontsize{14}{16pt}\selectfont $b$}, {\fontsize{14}{16pt}\selectfont $c$} is irrational,

\item[3.] One of the face diagonals {\fontsize{14}{16pt}\selectfont $d_{ab}$}, {\fontsize{14}{16pt}\selectfont $d_{bc}$}, {\fontsize{14}{16pt}\selectfont $d_{ac}$} is irrational.
\end{itemize}

In this paper we consider the third problem and show how parametric formulae for such NPC can be constructed.

PC problem is equivalent to the existence of an integer (rational) solution of the Diophantine system:
\vskip+0.02cm
{\fontsize{14}{16pt}\selectfont
\begin{align*}
    a^2+b^2 & =d_{ab}^{\,2}, \\[0.2cm]
    b^2+c^2 & =d_{bc}^{\,2}, \\[0.2cm]
    a^2+c^2 & =d_{ac}^{\,2}, \\[0.2cm]
    a^2+b^2+c^2 & =d_s^{\,2}.
\end{align*}
}
\vskip+0.2cm

In {\fontsize{14}{16pt}\selectfont [1]} are found complete parametrizations for NPC (only one face diagonal {\fontsize{14}{16pt}\selectfont $d_{ab}$} is irrational) by using a pair of solutions of the congruent number equation. Given parametrizations are not rational and depend on two parameters.

In {\fontsize{14}{16pt}\selectfont [2--3]} are summarized the literature written about PC. In numerous papers, only a few are related to the rational parametrizations {\fontsize{14}{16pt}\selectfont [4--9]}.

In this paper we use the same notations as used in {\fontsize{14}{16pt}\selectfont [1]}.

\vskip+1cm
\section{Basic Algebra}
\vskip+0.5cm

Theorem 3 in {\fontsize{14}{16pt}\selectfont [1]} establishes PC equivalent equation:
\vskip+0.02cm
{\fontsize{14}{16pt}\selectfont
\begin{equation}\label{eq:1}
    \Big(\frac{1-\gamma^2}{2\gamma}\Big)^2+\Big(\frac{1-\beta^2}{2\beta}\Big)^2=\Big(\frac{1-\alpha^2}{2\alpha}\Big)^2,
\end{equation}
}
\vskip+0.2cm
\noindent
where {\fontsize{14}{16pt}\selectfont $\alpha$}, {\fontsize{14}{16pt}\selectfont $\beta$}, {\fontsize{14}{16pt}\selectfont $\gamma$} are nontrivial rational numbers, i.e. from set
\vskip+0.02cm
{\fontsize{14}{16pt}\selectfont
    $$    \mathbb{Q}\setminus\{0;\pm\,1\}.      $$
    }
\vskip+0.2cm
\noindent
After simplification:
\vskip+0.02cm
{\fontsize{14}{16pt}\selectfont
    $$    \Big(\frac{1-\gamma^2}{2\gamma}\Big)^2=\frac{\big(1-(\alpha\beta)^2\big)\Big(1-\Big(\dfrac{\alpha}{\beta}\Big)^2\Big)}
                {4(\alpha\beta)\cdot\Big(\dfrac{\alpha}{\beta}\Big)}\,.     $$
    }
\vskip+0.2cm
\noindent
By using new notations:
\vskip+0.02cm
{\fontsize{14}{16pt}\selectfont
\begin{equation}\label{eq:2}
    \xi\equiv\alpha\beta, \quad \zeta\equiv\frac{\alpha}{\beta}\,;
\end{equation}
}
\vskip+0.2cm
\noindent
PC equation takes the form
\vskip+0.02cm
{\fontsize{14}{16pt}\selectfont
\begin{equation}\label{eq:3}
    \Big(\frac{1-\gamma^2}{2\gamma}\Big)^2=\frac{(1-\xi^2)(1-\zeta^2)}{4\xi\zeta}\,.
\end{equation}
}
\vskip+0.2cm
\noindent
It is clear that the product (ratio) {\fontsize{14}{16pt}\selectfont $\xi\zeta$} is square {\fontsize{14}{16pt}\selectfont $(=\square)$} and
\vskip+0.02cm
{\fontsize{14}{16pt}\selectfont
    $$    \xi,\,\zeta\in\mathbb{Q}\setminus\{0;\pm\,1\}, \quad \xi\neq\zeta.    $$
    }
\vskip+0.2cm
\noindent
From \eqref{eq:2} and \eqref{eq:3} we obtain:

\bigskip
\begin{theorem}\label{th:1}
The existence of Perfect cuboid is equivalent to the existence of nontrivial different {\fontsize{14}{16pt}\selectfont $\xi$} and {\fontsize{14}{16pt}\selectfont $\zeta$} rational numbers  satisfying three conditions:
\vskip+0.02cm
{\fontsize{14}{16pt}\selectfont
\begin{align}
    \xi\zeta & =\square, \label{eq:4} \\[0.2cm]
    (1-\xi^2)(1-\zeta^2) & =\square, \label{eq:5} \\[0.2cm]
    (1-\xi^2)(1-\zeta^2)+4\xi\zeta & =\square. \label{eq:6}
\end{align}
    }
\end{theorem}
\bigskip

Let's discuss the identity like condition \eqref{eq:5}:
\vskip+0.02cm
{\fontsize{14}{16pt}\selectfont
\begin{equation}\label{eq:7}
    (1-T^2)\big(1-(4T^3-3T)^2\big)=\big[(1-T^2)(1-4T^2)\big]^2.
\end{equation}
}
\vskip+0.2cm
\noindent
We try to replace {\fontsize{14}{16pt}\selectfont $\xi$} and {\fontsize{14}{16pt}\selectfont $\zeta$} with {\fontsize{14}{16pt}\selectfont $T$} and {\fontsize{14}{16pt}\selectfont $4T^3-3T$}, respectively.

From condition \eqref{eq:4}, it is possible if
\vskip+0.02cm
{\fontsize{14}{16pt}\selectfont
\begin{equation}\label{eq:8}
    4T^2-3=\square\,.
\end{equation}
}
\vskip+0.2cm
\noindent
The complete rational parametrization of \eqref{eq:8} is
\vskip+0.02cm
{\fontsize{14}{16pt}\selectfont
    $$  T=\frac{t^2+3}{4t}\,,       $$
}
\vskip+0.2cm
\noindent
where {\fontsize{14}{16pt}\selectfont $t$} is an arbitrary (nonzero) rational number.

We obtain the rational parametrization of conditions \eqref{eq:4} and \eqref{eq:5} by formulae:
\vskip+0.02cm
{\fontsize{14}{16pt}\selectfont
\begin{equation}\label{eq:9}
    \xi=\frac{t^2+3}{4t}\,, \quad \zeta=\frac{t^2+3}{4t}\,\Big(\frac{t^2-3}{2t}\Big)^2.
\end{equation}
}
\vskip+0.2cm
\noindent
Of course \eqref{eq:9} represents the incomplete parametrization of \eqref{eq:4} and \eqref{eq:5}, because identity \eqref{eq:7} does not describe all rational solutions of \eqref{eq:5}.

\vskip+1cm
\section{The First Parametrization}
\vskip+0.5cm

PC equation \eqref{eq:1} is obtained from the system {\fontsize{14}{16pt}\selectfont [1]}:
\vskip+0.02cm
{\fontsize{14}{16pt}\selectfont
\begin{equation}\label{eq:10}
\begin{gathered}
    \frac{d_s}{a}=\frac{1+\alpha^2}{2\alpha}\,, \quad \frac{d_{bc}}{a}=\frac{1-\alpha^2}{2\alpha}\,, \\[0.2cm]
    \frac{d_{ac}}{a}=\frac{1+\beta^2}{2\beta}\,, \quad \frac{c}{a}=\frac{1-\beta^2}{2\beta}\,, \\[0.2cm]
    \frac{d_{ab}}{a}=\frac{1+\gamma^2}{2\gamma}\,, \quad \frac{b}{a}=\frac{1-\gamma^2}{2\gamma}\,.
\end{gathered}
\end{equation}
}
\vskip+0.2cm
\noindent
From \eqref{eq:2} follows
\vskip+0.02cm
{\fontsize{14}{16pt}\selectfont
\begin{equation}\label{eq:11}
    \alpha^2=\xi\zeta, \quad \beta^2=\frac{\xi}{\zeta}\,,
\end{equation}
}
\vskip+0.2cm
\noindent
and if \eqref{eq:4} and \eqref{eq:5} are satisfied, then the ratios
\vskip+0.02cm
{\fontsize{14}{16pt}\selectfont
    $$  \frac{d_s}{a}\,, \;\; \frac{d_{bc}}{a}\,, \;\; \frac{d_{ac}}{a}\,, \;\; \frac{c}{a}\,, \;\; \frac{b}{a}     $$
}
\vskip+0.2cm
\noindent
are rationals. So,

\bigskip
\begin{theorem}\label{th:2}
Nearly-perfect cuboid $($only one face diagonal is irrational$)$ is obtained by nontrivial different {\fontsize{14}{16pt}\selectfont $\xi$} and {\fontsize{14}{16pt}\selectfont $\zeta$} rational numbers satisfying two conditions:
{\fontsize{14}{16pt}\selectfont
\begin{align*}
    \xi\zeta & =\square, \\[0.2cm]
    (1-\xi^2)(1-\zeta^2) & =\square.
\end{align*}
    }
\end{theorem}
\bigskip

From \eqref{eq:9} and \eqref{eq:11}:
\vskip+0.02cm
{\fontsize{14}{16pt}\selectfont
\begin{equation}\label{eq:12}
    \alpha=\frac{t^4-9}{8t^2}\,, \quad \beta=\frac{2t}{t^2-3}\,.
\end{equation}
}
\vskip+0.2cm
\noindent
By substituting \eqref{eq:12} in \eqref{eq:10} we get the first parametrization for NPC:
\vskip+0.02cm
{\fontsize{14}{16pt}\selectfont
\begin{align*}
    &\qquad \textbf{I\; parametrization} \\[0.5cm]
    a & =16t^2(t^4-9), \\[0.2cm]
    b & =(t^4-10t^2+9)(t^4+2t^2+9), \\[0.2cm]
    c & =4t(t^2+3)(t^4-10t^2+9); \\[0.5cm]
    d_{ac} & =4t(t^2+3)(t^4-2t^2+9), \\[0.2cm]
    d_{bc} & =(t^4-1)(t^4-81), \\[0.2cm]
    d_s & =t^8+46t^4+81.
\end{align*}
    }

\newpage
\section{The Second Parametrization}
\vskip+0.5cm

In {\fontsize{14}{16pt}\selectfont [1]} three PC equations are obtained so, naturally we expect other parametrizations.

Consider the PC equation from Theorem 1 in {\fontsize{14}{16pt}\selectfont [1]}:
\vskip+0.2cm
{\fontsize{14}{16pt}\selectfont
\begin{equation}\label{eq:13}
    \Big(\frac{2\alpha}{1+\alpha^2}\Big)^2+\Big(\frac{2\gamma}{1+\gamma^2}\Big)^2=\Big(\frac{2\beta}{1+\beta^2}\Big)^2.
\end{equation}
}
\vskip+0.2cm
\noindent
After simplification:
\vskip+0.2cm
{\fontsize{14}{16pt}\selectfont
    $$    \Big(\frac{2\gamma}{1+\gamma^2}\Big)^2=\frac{4\beta^2\big(1-(\alpha\beta)^2\big)\Big(1-\Big(\dfrac{\alpha}{\beta}\Big)^2\Big)}
                            {(1+\alpha^2)^2(1+\beta^2)^2}\,.        $$
}
\vskip+0.35cm
\noindent
PC equation \eqref{eq:13} is obtained from the following system {\fontsize{14}{16pt}\selectfont [1]}:
\vskip+0.2cm
{\fontsize{14}{16pt}\selectfont
\begin{equation}\label{eq:14}
\begin{gathered}
    \frac{a}{d_s}=\frac{2\alpha}{1+\alpha^2}\,, \quad \frac{d_{bc}}{d_s}=\frac{1-\alpha^2}{1+\alpha^2}\,, \\[0.35cm]
    \frac{b}{d_s}=\frac{1-\beta^2}{1+\beta^2}\,, \quad \frac{d_{ac}}{d_s}=\frac{2\beta}{1+\beta^2}\,, \\[0.35cm]
    \frac{c}{d_s}=\frac{2\gamma}{1+\gamma^2}\,, \quad \frac{d_{ab}}{d_s}=\frac{1-\gamma^2}{1+\gamma^2}\,.
\end{gathered}
\end{equation}
}
\vskip+0.2cm
\noindent
By substituting \eqref{eq:12} in \eqref{eq:14}, we obtain the second parametrization for NPC:
\vskip+0.02cm
{\fontsize{14}{16pt}\selectfont
\begin{align*}
     &\qquad \textbf{II\; parametrization} \\[0.5cm]
    a & =16t^2(t^4-9)(t^4-2t^2+9), \\[0.2cm]
    b & =(t^4-10t^2+9)(t^8+46t^4+81), \\[0.2cm]
    c & =4t(t^2-3)(t^4-10t^2+9)(t^4+2t^2+9); \\[0.5cm]
    d_{ac} & =4t(t^2-3)(t^8+46t^4+81), \\[0.2cm]
    d_{bc} & =(t^4-2t^2+9)(t^8-82t^4+81), \\[0.2cm]
    d_s & =(t^4-2t^2+9)(t^8+46t^4+81).
\end{align*}
    }

\vskip+1cm
\section{The Third Parametrization}
\vskip+0.5cm

From Theorem 2 in {\fontsize{14}{16pt}\selectfont [1]} PC equation is
\vskip+0.02cm
{\fontsize{14}{16pt}\selectfont
\begin{equation}\label{eq:15}
    \Big(\frac{2\gamma}{1-\gamma^2}\Big)^2+\Big(\frac{2\beta}{1-\beta^2}\Big)^2=\Big(\frac{2\alpha}{1-\alpha^2}\Big)^2.
\end{equation}
}
\vskip+0.2cm
\noindent
After simplification:
\vskip+0.02cm
{\fontsize{14}{16pt}\selectfont
    $$    \Big(\frac{2\gamma}{1-\gamma^2}\Big)^2=\frac{4\alpha^2\big(1-(\alpha\beta)^2\big)\Big(1-\Big(\dfrac{\beta}{\alpha}\Big)^2\Big)}
                    {(1-\alpha^2)^2(1-\beta^2)^2}\,.        $$
}
\vskip+0.2cm
\noindent
In this case insert
\vskip+0.02cm
{\fontsize{14}{16pt}\selectfont
    $$    \alpha\beta=\frac{t^2+3}{4t}\,, \quad \frac{\beta}{\alpha}=\frac{t^2+3}{4t}\,\Big(\frac{t^2-3}{2t}\Big)^2     $$
}
\vskip+0.2cm
\noindent
into the generating system {\fontsize{14}{16pt}\selectfont [1]} for PC equation \eqref{eq:15}:
\vskip+0.02cm
{\fontsize{14}{16pt}\selectfont
\begin{gather*}
    \frac{d_s}{a}=\frac{1+\alpha^2}{1-\alpha^2}\,, \quad \frac{d_{bc}}{a}=\frac{2\alpha}{1-\alpha^2}\,, \\[0.2cm]
    \frac{d_{ac}}{a}=\frac{1+\beta^2}{1-\beta^2}\,, \quad \frac{c}{a}=\frac{2\beta}{1-\beta^2}\,, \\[0.2cm]
    \frac{d_{ab}}{a}=\frac{1+\gamma^2}{1-\gamma^2}\,, \quad \frac{b}{a}=\frac{2\gamma}{1-\gamma^2}\,.
\end{gather*}
}
\vskip+0.2cm
\noindent
We receive the third parametrization for NPC:
\vskip+0.02cm
{\fontsize{14}{16pt}\selectfont
\begin{align*}
    &\qquad \textbf{III\; parametrization} \\[0.5cm]
    a & =(t^4-1)(t^4-81), \\[0.2cm]
    b & =4t(t^2-3)(t^4+2t^2+9), \\[0.2cm]
    c & =16t^2(t^4-9); \\[0.5cm]
    d_{ac} & =t^8+46t^4+81, \\[0.2cm]
    d_{bc} & =4t(t^2-3)(t^4+10t^2+9), \\[0.2cm]
    d_s & =(t^4-2t^2+9)(t^4+10t^2+9).
\end{align*}
    }

\newpage
\section{The Basis of a Computer Search}
\vskip+0.5cm

The goal of this paper is to find a perfect cuboid. The obtained parametrizations can form the basis for a computer search. The first, second and third parametrizations give for diagonal\; {\fontsize{14}{16pt}\selectfont $d_{ab}$}\; conditions of rationality, respectively:
\vskip+0.02cm
{\fontsize{14}{16pt}\selectfont
\begin{align*}
    \big[16t^2(t^4-9)\big]^2+\Big[(t^4-10t^2+9)(t^4+2t^2+9)\Big]^2 & =\square, \\[0.4cm]
    \Big[16t^2(t^4-9)(t^4-2t^2+9)\Big]^2+ \qquad\qquad\qquad & \\[0.2cm]
    +\Big[(t^4-10t^2+9)(t^8+46t^4+81)\Big]^2 & =\square, \\[0.4cm]
    \big[(t^4-1)(t^4-81)\big]^2+\Big[4t(t^2-3)(t^4+2t^2+9)\Big]^2 & =\square.
\end{align*}
}
\vskip+0.2cm
\noindent
If, by using a computer, it is possible to find the nontrivial 
$(\neq 0,\pm1,\pm3)$ rational {\fontsize{14}{16pt}\selectfont $t$}, for which one of the given expressions is square, then PC problem has a solution. Otherwise from Theorem \ref{th:1} it must be proved that two nontrivial different rational numbers which satisfied the \eqref{eq:4}--\eqref{eq:6} conditions do not exist.

\vskip+1cm
\begin{center}
     {\fontsize{14}{16pt}\selectfont \textbf{References} }
\end{center}
\vskip+0.5cm

\begin{enumerate}
\item[1.] Meskhishvili M., Perfect cuboid and congruent number equation solutions. \texttt{http://arxiv.org/pdf/1211.6548v2.pdf}

\medskip

\item[2.] van Luijk R., On perfect cuboids. \emph{Doctoraalscriptie, Universiteit \linebreak   Utrecht}, 2000. \texttt{http://www.math.leidenuniv.nl/reports/ \\ 2001-12.shtml}

\medskip

\item[3.] Meskhishvili M., Three-Century problem. \emph{Tbilisi}, 2013.

\medskip

\item[4.] Leech J., The rational cuboid revisited. \emph{Amer. Math. Monthly}  \textbf{84}  (1977), No. 7, 518--533.

\medskip

\item[5.] MacLeod A. J., Parametric expressions for a ``nearly-perfect'' cuboid. \texttt{http://ru.scribd.com/doc/58729272/Parametric-Equations- \\ for-Nearly-Perfect-Cubiods}
\medskip

\item[6.] Bromhead T., On square sums of squares. \emph{Math. Gaz.}  \textbf{44}  (1960), 219--220.

\medskip

\item[7.] Colman W. J. A., On certain semiperfect cuboids. \emph{Fibonacci Quart.}  \textbf{26}  (1988), No. 1, 54--57.

\medskip

\item[8.] Narumiya N., Shiga H., On certain rational cuboid problems. \emph{Nihonkai Math. J.}  \textbf{12}  (2001), No. 1, 75--88.

\medskip

\item[9.] Top J., Yui N., Congruent number problems and their variants.  \emph{Al\-go\-ri\-thmic number theory: lattices, number fields, curves and cryp\-tog\-ra\-phy},  613--639, \emph{Math. Sci. Res. Inst. Publ.}, 44, \emph{Cambridge Univ. Press, Cambridge}, 2008.

\end{enumerate}

\vskip+1.5cm

\noindent Author's address:

\medskip

\noindent {Georgian-American High School, 18 Chkondideli Str., Tbilisi 0180, Georgia.}

\noindent {\small \textit{E-mail:} \texttt{director@gahs.edu.ge} }

\end{document}